\documentclass[a4paper]{article}

\usepackage[T1]{fontenc}
\usepackage[utf8]{inputenc}

\usepackage{amsmath}
\usepackage{mathtools}
\usepackage{amssymb}

\usepackage{graphicx}  

\DeclarePairedDelimiter\floor{\lfloor}{\rfloor}

\begin{document}

\title{A hybrid Fourier-Prony method}
\author{Matteo Briani, Annie Cuyt, Wen-shin Lee}
\date{}
\maketitle

\begin{abstract}
  The FFT algorithm that implements the discrete Fourier transform is considered one of the top ten algorithms of the $20$th century.
Its main strengths are the low computational cost of $\mathcal{O}(n \log n$) and its stability.
It is one of the most commonly used algorithms to analyze signals with a dense frequency representation.
In recent years there has been an increasing interest in sparse signal representations and a need for algorithms that exploit such structure.
We propose a new technique that combines the properties of the discrete Fourier transform with the sparsity of the signal.
This is achieved by integrating ideas of Prony's method into Fourier's method.
The resulting technique has the same frequency resolution as the original FFT algorithm but uses fewer samples and can achieve a lower computational cost.
Moreover, the proposed algorithm is well suited for a parallel implementation. 
\end{abstract}

\small{\textbf{Keywords:} Fourier transform, Prony's method, sparse representation}

\small{\textbf{AMS-classification number:} 65T50}

%----------------------------------------------------%
% INTRODUCTION
%----------------------------------------------------%
\section{Introduction}\label{intro}

The Discrete Fourier Transform (DFT) is a mathematical transform that has been widely used from the beginning of the digital era.
It is applied in almost all fields of digital signal processing.
In the last century a big effort was to develop fast implementations. 
It is now possible to analyze large streams of samples using procedures of low computational complexity.
Many of the publications in this area are dealing with signals that have a dense structure in the frequency domain. 
Now, at the new frontier, one is considering the case of signals that have a sparse frequency representation.
Hence another step forward is possible if sparsity is exploited.
A recent review of such methods is presented in \cite{gilbert_recent_2014}.

One of these approaches involves the use of so-called parametric methods, such as \cite{schmidt_multiple_1986,roy_esprit_1989,hua_matrix_1990}.
The frequency resolution of these methods is not bounded, as in the discrete Fourier method, by the amount and the time span of the available samples. 
Thus they seem to be the first choice for exploiting sparse signals.
However, some extra care has to be taken due to their sensitivity to noise \cite{batenkov_stability_2016}.
For instance, in \cite{heider_sparse_2013-1} the authors use a Prony method in conjunction with an appropriate filtering technique to replace a standard DFT.
In particular, the original signal is first filtered into several frequency bands that contain only a portion of the signal components.
Afterwards, these components are moved apart in the spectrum by means of a permutation of the filtered samples.
A parametric method such as MUSIC \cite{schmidt_multiple_1986}, ESPRIT \cite{roy_esprit_1989} or Matrix Pencil \cite{hua_matrix_1990} is then used to detect the meaningful components.

Another approach falls under the sub-Nyquist sampling techniques.
In \cite{christlieb_multiscale_2015}, for the computation of the DFT, the signal is sampled at a rate that does not obey the Shannon-Nyquist (S-N) theorem \cite{shannon_communication_1949}.
This permits to collect fewer samples for the analysis but has the downside of causing aliasing, i.e.\ remapping of the signal components in a lower part of the spectrum.
Once the meaningful components are retrieved, they are reallocated in the proper positions thanks to the information extracted from an additional set of collected samples.
This second set is close enough to the first one to satisfy the S-N theorem and thus it allows to resolve the aliasing.
Of the same kind is the use of the Chinese Reminder Theorem (CRT) in employing sub-Nyquist sets of samples collected with different undersampling rates \cite{iwen_combinatorial_2010}.
If the undersampling rates are kept coprime, it is possible to resolve the aliasing by means of the CRT and obtain an alias-free Fourier transform.

In \cite{hassanieh_simple_2012} a probabilistic approach is used.
The signal samples are (pseudo)-randomly permuted and then filtered.
A short version of the DFT is computed and the biggest peaks are selected and their locations stored.
By repeating these steps with different pseudo-random permutations, it is possible to detect the true frequency locations, using a probabilistic argument.
Along the same line, are the papers \cite{hassanieh_nearly_2012, indyk_nearly_2014}.

In our method we combine different approaches, namely Fourier and Prony related methods.
We use Fourier techniques to process undersampled signals affected by aliasing.
In a second stage, small Prony systems are built in order to resolve the introduced aliasing and identify the non-aliased solution.
The parametric methods are not used to directly compute the aliased-free components, instead as a tool to extract information.
We make good use of the Chinese Remainder Theorem, but unlike in other works, we do not need to collect our samples at different sampling rates.
So in all, we develop a brand new method that processes sub-Nyquist samples with almost no loss in the frequency resolution.
Recently in \cite{potts_efficient_2016}, a related approach was discussed: there the combination of Fourier and Prony techniques relies on the fact that several Fourier transforms are computed with an offset of a single sample.
The current independently developed paper presents a more general setting where both Fourier and Prony can be affected by aliasing and, in particular, the Fourier transforms can be computed with a bigger offset allowing less restrictions on the data acquisition process.

The paper is organized as follows: Section \ref{under_time_shifts} deals with the effect of undersampling and treats time-shifted signals.
In Section \ref{method} we present the core of our technique and in Section \ref{colliding} we address the problem of discerning frequencies colliding because of possible aliasing.
In Section \ref{noise} we explain how the method copes with noise and, in Section \ref{numerical} we test our method on a numerical example.

%----------------------------------------------------%
% UNDERSAMPLING AND TIME SHIFTS
%----------------------------------------------------%
\section{Undersampling and time shifts}\label{under_time_shifts}

Let us consider a finite set of samples collected on a uniform time grid from a function $ \mathcal{X}: \left[a, b \right] \rightarrow \mathbb{C}$.
The cardinality of the set is $N$ and the acquired samples are indicated by $\boldsymbol{x} = (x_l)_{l=0}^{N-1}$.
The DFT coefficients $\boldsymbol{X} = (X_j)_{j=0}^{N-1}$ are defined by

\begin{equation}
  X_j := \sum_{l=0}^{N-1} x_l \mathrm{e}^{\frac{-2 \pi \mathrm{i}}{N} l j}, \qquad \mathrm{i}^2 = -1.
\end{equation}

Each Fourier coefficient $X_j$ is the inner product $\langle \boldsymbol{x}, \exp({\frac{-2 \pi \mathrm{i}}{N} j \ell })_{\ell=0}^{N-1} \rangle$. 
A Fourier coefficient $X_j$ can be regarded as the amount of a specific complex exponential that is present in the discrete signal $\boldsymbol{x}$.
We say that each $X_j$ is associated to a frequency, meaning the associated complex exponential.

The set of samples $\boldsymbol{x}$ can be computed back from the Fourier coefficients $\boldsymbol{X}$ using the Inverse Discrete Fourier Transfrom (IDFT) by

\begin{equation}
  x_l := \frac{1}{N} \sum_{j=0}^{N-1} X_j \mathrm{e}^{\frac{2 \pi \mathrm{i}}{N} l j}.
  \label{IDFT}
\end{equation}

The Shannon-Nyquist theorem \cite{shannon_communication_1949} states that we are able to exactly reconstruct bandlimited signals that have the same bandwidth as the sampling rate.
If the analyzed signal contains frequencies that are above the highest frequency we can retrieve, we encounter aliasing. 
We take care of this problem below. 
If not stated otherwise, we assume that the signal $\boldsymbol{x}$ satisfies the S-N theorem. 

Before proceeding, we also assume that the signal $\boldsymbol{x}$ has a sparse reprentation in the frequency domain, i.e. it has only $K\ll N$ non-vanishing Fourier coefficients $X_{j_k} \in \mathbb{C}, j_k \in \{0,\ldots,N-1\}, k=1,\ldots,K$.
In case the signal is perturbed by noise, we set a threshold $T \in \mathbb{R}$ and we assume that only $K$ Fourier coefficients have an amplitude larger than $T$, 
\begin{equation}
\lvert X_{j_k} \rvert \geq T.
\label{thresholdT}
\end{equation}

We now consider the following undersampled (possibly sub-Nyquist) version of the signal $\boldsymbol{x}$

\begin{equation}
  {_u}\boldsymbol{x} :=(x_{ul})_{l=0}^{\floor{N-1/u}},
\end{equation}

with $u \in \mathbb{N}$ and $\floor{}$ indicating the $\operatorname{floor}$ function.
We denote this smaller set of samples as $ {_u}\boldsymbol{x}$ and by ${_u}\boldsymbol{X}$ its associated DFT.
The frequency resolution of the DFT is given by the sampling rate over the number of analyzed samples.
A wider distance between samples implies that the highest frequency we are able to retrieve is now smaller and some coefficients may have been remapped to other frequencies.
Given a specific Fourier coefficient, there is no way to know if it has been affected by aliasing or not. 

Let us consider a Fourier coefficient $X_{j_k}$ and compute the undersampled DFT ${_u}\boldsymbol{X}$. 
The same coefficient $X_{j_k}$, due to aliasing, may now appear at a different index $\tilde{j_k}$.
In particular 

\begin{equation}
  X_{j_k} = {_u}X_{\tilde{{j_k}}} \quad  j_k = { u \tilde{{j_k}}}  \operatorname{mod}(N).
  \label{aliasingmod}
\end{equation}

However, in order to leave the notation lighter, we do not explicitly put a tilde over the aliased indices.
Moreover, when we are dealing with aliasing, it can happen that a non-vanishing Fourier coefficient is mapped down to an exponential with another associated non-vanishing coefficient. 
In this case the retrieved coefficient ${_u}X_{j_k}$ is the sum of two or more coefficients of the full-length DFT.
We refer to this phenomenon as frequency collision.
For the moment we assume that no frequency collision occurs in our examples, we introduce it in a later part of the paper.

We denote a shifted version of the signal $\boldsymbol{x}$ by

\begin{equation}
  {^s}\boldsymbol{x} := (x_{l+s})_{l=0}^{N-1} \quad s \in \mathbb{N}.
\end{equation}

We recall that, since the signal is assumed to be periodic, $ x_{l+\lambda N} = x_l, \lambda \in \mathbb{N}$.
Assuming that $\boldsymbol{X}$ and ${^s}\boldsymbol{X}$ are the noise-free Fourier coefficients computed from the set of samples $\boldsymbol{x}$ and ${^s}\boldsymbol{x}$ respectively, we have

\begin{align}
  {^s}{X}_j  &= \sum_{l=0}^{N-1} {^s}x_{l} e^{\frac{-2 \pi \mathrm{i}}{N} l j}, \nonumber \\
             &= \sum_{l=-s}^{N-1+s} {^s}x_{l} e^{\frac{-2 \pi \mathrm{i}}{N} l j}, \nonumber \\
             &= e^{ \frac{2 \pi \mathrm{i}}{N} sj} \sum_{l=-s}^{N-1+s} {^s}x_{l} e^{\frac{-2 \pi \mathrm{i}}{N} l j}  e^{ -\frac{2 \pi \mathrm{i}}{N} s j}, \nonumber \\
             &= e^{ \frac{2 \pi \mathrm{i}}{N} sj} \sum_{\bar{l}=0}^{N-1} x_{\bar{l}} e^{\frac{-2 \pi \mathrm{i}}{N} \bar{l} j} \left( \text{with} \quad \bar{l} = l + s \right), \nonumber \\
             &= e^{ \frac{2 \pi \mathrm{i}}{N} s j} X_j. \label{time_shift_derivation}
\end{align}

In particular, the shifted set has the same Fourier coefficients $X_j$ multiplied by the complex exponential $ e^{\frac{2 \pi \mathrm{i}}{N} s j}$. 
In the time domain this corresponds to a phase variation of the signal.

From \eqref{aliasingmod} and \eqref{time_shift_derivation} we see that for a shifted undersampled version of the signal it holds

\begin{equation}
  {^s_u}X_j = e^{\frac{2 \pi \mathrm{i}}{N} s j} {_u}X_j.
  \label{time_shift_under}
\end{equation}

This is particularly important because, given an index $j$, each decimated DFT equals ${_u}\boldsymbol{X}$ times a phase component.
Moreover, in the case of a frequency collision, \eqref{time_shift_under} implies that the collision index remains the same for each decimated ${_u^s}\boldsymbol{X}$.
In Section \ref{method} we use \eqref{time_shift_under} to resolve the aliasing issue.
We underline that \eqref{time_shift_under} holds for any value of $s$.
If noise is added to the Fourier coefficients, the relation between ${^s}X_j$ and $X_j$ is not exact anymore due to the non-periodicity of the noise.

Let us consider the sets of samples $\boldsymbol{x}$ and ${^1}\boldsymbol{x}$ $(s=1)$, and compute their respective non-vanishing Fourier coefficients $X_{j_k}$ and ${^1}{X}_{j_k}$, $ k=1,\ldots,K$.
One implication of the previous consideration is that dividing each of the non-vanishing coefficients ${^1}{X}_{j_k}$ by $X_{j_k}$ we obtain $e^{\frac{2 \pi \mathrm{i}}{N} j_k}$.
In other words, given the two sets of Fourier coefficients $\boldsymbol{X}$ and ${^1}{\boldsymbol{X}}$ we can extract their associated complex exponentials.
If we consider the Fourier transform of $\boldsymbol{x}$ and ${^s}\boldsymbol{x}$, the division of ${^s}{X}_{j_k}$ by $X_{j_k}$ leads to $e^{\frac{2 \pi \mathrm{i}}{N} s j_k}$. 
Even if it does not seem useful at the moment, we make good use of this property in the sequel.
We underline that the information on the associated complex exponential is already hidden inside the Fourier coefficients ${^s}\boldsymbol{X}$.
However, this property holds only when the Nyquist rate is satisfied ( \cite{christlieb_multiscale_2015} and \cite{potts_efficient_2016} are following this approach).
If a Fourier coefficient is affected by aliasing, then we are left with a set of plausible solutions 

\begin{equation}
  S_{j_k} := \left\{\exp \left( \frac{2 \pi \mathrm{i}}{N} s j_k + \frac{2 \pi \mathrm{i}}{s} \ell \right)\right\}_{\ell=0}^{s-1}.
  \label{solutionsS}
\end{equation}

We call generator of the set $S_{j_k}$ the exponential $\exp( \frac{2 \pi \mathrm{i}}{N} s j_k)$.
Again, it is not possible to discern which is the right one.
However, using particular values of $u$ and $s$ solves this problem.

%----------------------------------------------------%
% FIXING THE ALIASING
%----------------------------------------------------%
\section{Fixing the aliasing}\label{method}

As stated in Section \ref{under_time_shifts}, each Fourier coefficient of ${_u}\boldsymbol{X}$ might have been affected by aliasing. 
For each index $j_k$ we define  the set

\begin{equation}
  U_{j_k} := \left\{ \exp \left( \frac{2 \pi \mathrm{i}}{N} u j_k + \frac{2 \pi \mathrm{i}}{u}\ell \right) \right\}_{\ell=0}^{u-1}
\end{equation}

and we call the exponential $\exp{\left( \frac{2 \pi \mathrm{i}}{N} u j_k \right)}$ its generator.

A smart way of using the parameters $u$ and $s$ can fix the aliasing.
The key point of this technique is to choose $u$ and $s$ to be coprime ( more general discussion in \cite{Cu.Le:ana:17}).
The amount of undersampling $u$ or the shift factor $s$ does not matter as long as this requirement is satisfied.

\begin{figure}
  \centering
  \includegraphics[scale=0.5]{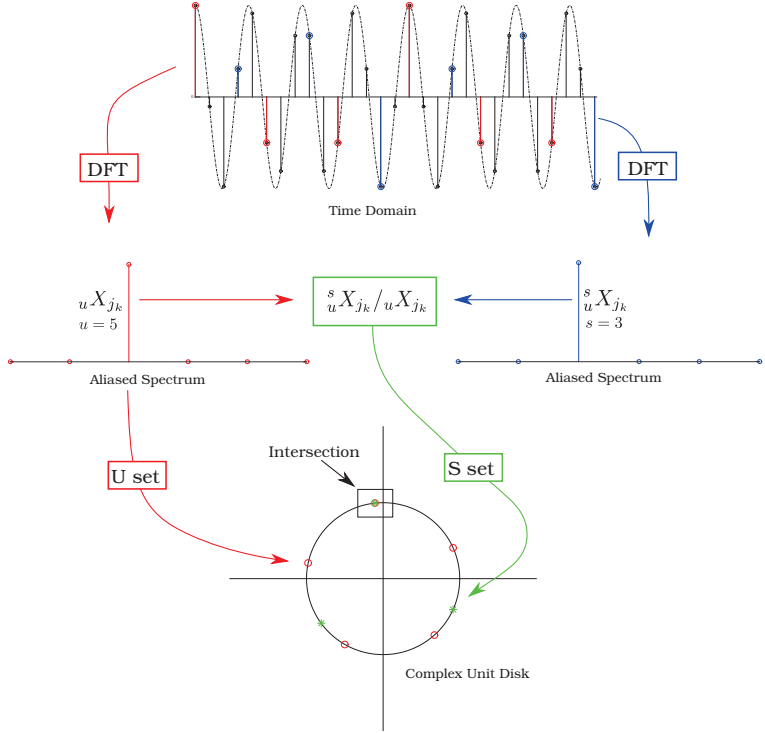}
  \caption{\it The core idea for solving the aliasing when a DFT is computed using a sub-Nyquist set of samples.
  However, if the Fourier coefficients are corrupted by noise, the division ${_u^s}X_{j_k}/{_u}X_{j_k}$ might not yield the right result.
  See Section \ref{colliding} for a discussion and solution of the problem.}
  \label{methodFigure}
\end{figure}

The procedure goes as follows: we calculate the DFT of ${_u}\boldsymbol{x}$ and we compute, for each non-vanishing coefficient, the set $U_{j_k} $.
We then consider the time-shifted set of samples ${_u^s}\boldsymbol{x}$ and compute the Fourier coefficients ${_u^s}\boldsymbol{X}$.
For each non-vanishing $j_k$ the division ${_u^s}X_{j_k}/{_u}X_{j_k}$ returns the exponential $\exp( \frac{2 \pi \mathrm{i}}{N} s j_k)$.
However, as previously stated,  $\exp( \frac{2 \pi \mathrm{i}}{N} s j_k)$ indicates the set of plausible exponentials \eqref{solutionsS}.
Each set contains the non-aliased exponential and, since $u$ and $s$ are coprime, these two sets share one and only one element, which is the right complex exponential associated to the Fourier coefficient $j_k$ (see also \cite{Cu.Le:how:17}).
Figure \ref{methodFigure} summarizes this technique.

For each non-vanishing $j_k$ we find the intersection between the two sets $U_{j_k}$ and $S_{j_k}$.
In order to do so, two ways are possible: search the intersection computing the distance among their elements, or use the Euclidean algorithm.
For the first approach we build a distance matrix from the two sets and simply extract the exponentials that are closer together. 
However, this might not be the best approach. 
Indeed, we have to compute all the distances among the possible solutions and then search for the smallest one.
So the computational cost increases for larger values of $u$ and $s$.

A less computationally demanding approach is based on the Euclidean algorithm. 
We recall that, since $u$ and $s$ are coprime, the Euclidean algorithm states that there exist two integers $t$ and $v$ that satisfy the equation $ut + sv = 1$.
The integer pair $t$, $v$ is not unique and might be chosen accordingly to the noise level as we see later.
We then consider the generators of the sets $U_{j_k}$ and $S_{j_k}$; namely $\exp( \frac{2 \pi \mathrm{i}}{N} u j_k)$ and $\exp( \frac{2 \pi \mathrm{i}}{N} s j_k)$.
We multiply each exponent by, respectively, $t$ and $v$ obtaining 

\begin{equation}
  e^{\frac{2 \pi \mathrm{i}}{N} (ut+sv) j_k} = e^{\frac{2 \pi \mathrm{i}}{N} j_k} 
\end{equation}

which is the solution we are seeking.
However, the downside of this approach is the sensitivity on the noise.
In fact, when we are computing $e^{\frac{2 \pi \mathrm{i}}{N} (ut+sv) j_k}$ we are powering $\exp( \frac{2 \pi \mathrm{i}}{N} u j_k)$ and $\exp( \frac{2 \pi \mathrm{i}}{N} s j_k)$.
Since $\exp( \frac{2 \pi \mathrm{i}}{N} u j_k)$ and $\exp( \frac{2 \pi \mathrm{i}}{N} s j_k)$ are computed from ${_u^s}X_{j_k}/{_u}X_{j_k}$, they are corrupted by the noise present in the Fourier coefficients.
Therefore, when powering the generators of $U_{j_k}$ and $S_{j_k}$, we are also powering the noise, expect when we can find a suitable pair of integers $t$, $v$ with small absolute values.

We also point out that the second set of plausible solutions $S_j$ can be computed without performing an additional Fourier transform.
We do not need to compute ${_u^s}\boldsymbol{X}$ but just solve a linear system of equations involving the non-vanishing ${_u}X_{j_k}$.
In order to do so we build a Vandermonde matrix $V$ from the $K$ non-vanishing ${_u}X_{j_k}$, 

\begin{equation}
  V = 
  \begin{pmatrix}
  1                   &  1                   &  1                   &  \cdots  &  1                   \\
  {_u}X_{j_{k_1}}  &  {_u}X_{j_{k_2}}  &  {_u}X_{j_{k_3}}  &  \cdots  &  {_u}X_{j_{k_K}}  \\
  {_u}X_{j_{k_1}}^2  &  {_u}X_{j_{k_2}}^2 &  {_u}X_{j_{k_3}}^2    &  \cdots  &  {_u}X_{j_{k_K}}^2   \\
    \vdots & & & & \vdots \\
    {_u}X_{j_{k_1}}^{K-1}  &  {_u}X_{j_{k_2}}^{K-1} &  {_u}X_{j_{k_3}}^{K-1}    &  \cdots  &  {_u}X_{j_{k_K}}^{K-1}  \\
  \end{pmatrix} 
\end{equation}

and consider the vector $ \tilde{\boldsymbol{x}} := (x_{ul+s})_{l=0}^{K-1}$.
We seek the solution $\boldsymbol{y}$ of the linear system $V \boldsymbol{y} = \tilde{\boldsymbol{x}}$.
The vector $\boldsymbol{y}$ is of the form ${_u^s}{X}_{j_k}$, $k = 1,\ldots,K$.
However, the Vandermonde matrix $V$ may be ill-conditioned, due to the $K$ selected Fourier coefficients on the complex unit circle.
How this can be overcome is indicated in \cite{Cu.Le:spa:16}.

We proceed with an example to better illustrate the phenomenon.
Suppose we have $N = 1000$ samples of the discretized signal $\boldsymbol{x}$. 
Let's assume it has $4$ meaningful Fourier coefficients $\exp \left(\frac{2 \pi \mathrm{i}}{N} j_k \right)$, with $j_k = 11, 22, 33, 44$.
The condition number of the associated Vandermonde matrix is bigger than $10^4$.
On the other hand, if the undersampling rate $u$ is chosen to be $250$, the Fourier coefficients are associated with the complex exponentials $\exp{ \left(\frac{2 \pi \mathrm{i}}{N} j_k u \right)} = \exp{ \left(\frac{2 \pi \mathrm{i}}{4} \tilde{j_k} \right)}, \tilde{j_k} = 1, 2, 3, 0$, which are equally spaced on the complex unit circle and form a Vandermonde matrix $V$ with condition number $1$.
As we can see, choosing an appropriate undersampling rate may lead to a better conditioning of the matrix $V$ thus making this approach more effective and even computationally less expensive than computing the shifted Fourier transform.

%----------------------------------------------------%
% COLLIDING FREQUENCIES 
%----------------------------------------------------%
\section{Colliding frequencies}\label{colliding}

When analyzing an undersampled version of a discrete signal $\boldsymbol{x}$, we may encounter aliasing and some Fourier coefficients may be mapped to the wrong complex exponential.
When no frequency collision occurs, the proposed method guarantees the correct remapping of the Fourier coefficients affected by aliasing.
However this ideal situation rarely happens due to the Fourier leakage effect.

The leakage effect appears when, in the signal $\boldsymbol{x}$, a component cannot be exactly represented by a single complex exponential of the form $\exp( \frac{2 \pi \mathrm{i}}{N} j), j \in \mathbb{N}$.
In this case, the specific frequency leaks to neighbouring frequencies and it affects all Fourier coefficients \cite{thompson_leakage_1980}.
The effect is more evident from the Fourier coefficients near the location of the exact frequency.

The situation becomes problematic when we consider the undersampled signal ${_u}\boldsymbol{x}$.
Due to the aliasing phenomenon, the Fourier coefficients are mapped to other complex exponentials. 
Moreover, if the leakage effect is present, a frequency collision will occur because of the component that leaked over all coefficients.
This is limiting the applicability of the proposed method due to the inexactness result that appears when dividing ${_u^s}{X}_j$ by ${_u}X_j$.

Let $X_{j_1}, X_{j_2}$ be non-vanishing Fourier coefficients of $\boldsymbol{x}$.
We consider ${_u}\boldsymbol{X}$ and we assume the coefficients are colliding at the index $\hat{j}$, generating ${_u}X_{\hat{j}}$.
We also compute the DFT of ${_u^s}\boldsymbol{x}$ that returns the coefficient ${_u^s}X_{\hat{j}}$.
We recall that, from \eqref{time_shift_under}, the shift factor $s$ does not influence the index where the frequencies are colliding, it only affects the phases of each Fourier coefficient.
When dividing the two we obtain

\begin{align}
  \frac{_u^s{}{X}_{\hat{j}}}{ {_u}X_{\hat{j}}} &= \frac{ {^s}X_{j_1} + {^s}X_{j_2} }{ {X}_{j_1} + {X}_{j_2} } = \frac{ e^{\frac{2 \pi \mathrm{i}}{N} s j_1}X_{j_1} + e^{\frac{2 \pi \mathrm{i}}{N} s j_2}X_{j_2} }{ {X}_{j_1} + {X}_{j_2} }
\end{align}

and we are unable to extract the complex exponentials $e^{\frac{2 \pi \mathrm{i}}{N} s j_1}, e^{\frac{2 \pi \mathrm{i}}{N} s j_2}$.
Without the correct information about the exponentials $e^{\frac{2 \pi \mathrm{i}}{N} s j_1}, e^{\frac{2 \pi \mathrm{i}}{N} s j_2}$ we are not able to apply the proposed method and fix the aliasing.

However, it is still possible to make good use of ${_u^s}X_{\hat{j}}$.
We consider the following discrete signals ${_u}\boldsymbol{x}, {_u^s}\boldsymbol{x}, \ldots, {_u^{(M-1)s}}\boldsymbol{x}$ and their relative DFTs ${_u}\boldsymbol{X},{_u^s}\boldsymbol{X}, \ldots, {_u^{(M-1)s}}\boldsymbol{X}$.
It remains true that the frequencies $j_1, j_2$ are colliding, now for all $m$ at the index $\hat{j}$ of ${_u^{ms}}\boldsymbol{X}$.
We rewrite each Fourier coefficient ${_u^{ms}}X_{\hat{j}}$ as 

\begin{align}
  _u{X}_{\hat{j}} &= {{X}}_{j_1} + {{X}}_{j_2}, \nonumber  \\
  {_u^s{X}}_{\hat{j}} &= e^{ \frac{2 \pi \mathrm{i}}{N} s j_1 } {{X}}_{j_1} + e^{ \frac{2 \pi \mathrm{i}}{N} s j_2 } {{X}}_{j_2}, \nonumber  \\
                    &\vdots \nonumber\\
  {_u^{(M-1)s}{X}}_{\hat{j}} &= e^{ \frac{2 \pi \mathrm{i}}{N} (M-1) s j_1 } {{X}}_{j_1} + e^{ \frac{2 \pi \mathrm{i}}{N} (M-1) s j_2 } {{X}}_{j_2}. 
  \label{Pseq}
\end{align}

From \eqref{Pseq} we obtain the sequence

\begin{equation}
  P(m) :=  e^{ \frac{2 \pi \mathrm{i}}{N} m s j_1 } {{X}}_{j_1} + e^{ \frac{2 \pi \mathrm{i}}{N} m s j_2 } {{X}}_{j_2}
  \label{Psequence}
\end{equation}

From \eqref{Psequence} we aim to extract the exponentials $e^{\frac{2 \pi \mathrm{i}}{N} s j_1}, e^{\frac{2 \pi \mathrm{i}}{N} s j_2}$.
This is possible using a parametric method, in particular one of the Prony's family \cite{roy_esprit_1989,hua_matrix_1990,schmidt_multiple_1986}.
Prony's method and its variances are suited to analyze sums of complex exponentials.
More precisily, given the sequence $P(\cdot)$, the method returns the coefficients ${_u}X_{j_1}, {_u}X_{j_2}$ and their associated exponentials $e^{\frac{2 \pi \mathrm{i}}{N} s j_1}, e^{\frac{2 \pi \mathrm{i}}{N} s j_1}$.
Their frequency resolution is not restricted to a pre-assigned grid.
It is indeed mandatory to have a high frequency resolution in this part of the method.
Any loss of precision in the computation of the exponentials will lead to the wrong solution.
On the other hand, drawbacks of these methods are the computational cost and the susceptibility to high level of noise.  In Section \ref{noise} we show how to cope with a sequence $P(\cdot)$ corrupted by noise.

Given the exponentials $e^{\frac{2 \pi \mathrm{i}}{N} s j_1}, e^{\frac{2 \pi \mathrm{i}}{N} s j_2}$ and the coefficients ${_u}X_{j_1}, {_u}X_{j_2}$ we are able to resolve the aliasing for the coefficient $j_1, j_2$ and to split the coefficient ${_u^s}X_{\hat{j}}$ into its two components.
However, most of the time it is not necessary to do so.
If the collision was caused by the leakage of one of the frequencies, the absolute value of one of the coefficients $X_{j_1}, X_{j_2}$ may be small.
If $\lvert X_{j_1} \rvert$ is very small, the frequency $j_1$ is not close to the true leaking frequency.
Its contribution to another Fourier coefficient is bigger and its presence is not affecting the retrievial of the other frequency.
However, it is mandatory to have a criterion to discern wether a collision has happened or not.

In \cite{christlieb_multiscale_2015} the authors state a necessary condition for no collision to occur.
They do not use the sequence $P(\cdot)$ but only the original ${_u}X_j$ and the shifted Fourier bin ${_u^s}X_j$.
No collision occurs on bin $j$ only if 

\begin{equation}
  \label{nocollision}
  \lvert  \frac{ {_u^s}X_j }{ {_u}X_j } \rvert = \lvert e^{ \frac{2 \pi \mathrm{i}}{N} s j } \rvert = 1 .
\end{equation}

If a collision occurs, \eqref{nocollision} is no longer satisfied.
It is however more straightforward and reliable to analyze the sequence $P(\cdot)$.
In fact, from a Hankel matrix built from the sequence $P(\cdot)$ and its singular value decomposition, we can extract the number of meaningful components thus detecting a collision \cite{Cu.Ts.ea:fai:17}.
For an example see Section \ref{numerical}.

%----------------------------------------------------%
% NOISE
%----------------------------------------------------%
\section{Noise}\label{noise}

Before proceeding to the case of signals corrupted by noise, we recap the steps of the proposed method.

Given a signal $\boldsymbol{x}$ we fix an undersampling rate $u$, a shift factor $s$ coprime with $u$, and the number of shifted and undersampled DFT coefficients $M$.
We then consider ${_u}\boldsymbol{X}$, ${_u^s}\boldsymbol{X}$, \ldots, ${_u^{(M-1)s}}\boldsymbol{X}$.
The peaks of ${_u}\boldsymbol{X}$ correspond to a sum of one or more frequencies of $\boldsymbol{X}$ that can be affected by aliasing and collision.
Each peak indicates a set of plausible non-aliased frequencies $U_{\hat{j}}$. 
For fixed peak index $\hat{j}$ in ${_u}\boldsymbol{X}$, ${_u^s}\boldsymbol{X}$, \ldots, ${_u^{(M-1)s}}\boldsymbol{X}$, we form the sequence $P_{\hat{j}}(m)$ according to \eqref{Pseq} and \eqref{Psequence}.
Using Prony's method we analyze $P_{\hat{j}}(m)$ and extract its components.
For each component we are left with the set of plausible solutions $S_{j_1}, S_{j_2}, \ldots, S_{j_v}$, where $v$ is the number of detected components.
The intersection of each set $S_{j_1}, S_{j_2}, \ldots, S_{j_v}$ with the set $U_{\hat{j}}$ returns the true locations of the non-aliased frequencies.
Finally, the amplitude of each retrieved Fourier coefficient corresponds to the amplitude of the relative components in the sequence $P_{\hat{j}}(m)$.

At last, we consider the case when the signal $\boldsymbol{x}$ is corrupted by complex white gaussian noise $\boldsymbol{n} = (n)_{j=0}^{N-1}$.
We recall that the signal is still considered $K$-sparse if the threshold $T$ introduced in Section \ref{under_time_shifts} is such that \eqref{thresholdT} holds.

Because of the noise, the Fourier coefficients ${_u}\boldsymbol{X}$, ${_u^s}\boldsymbol{X}$, \ldots, ${_u^{(M-1)s}}\boldsymbol{X}$, are all corrupted.
This implies the corruption of the sequence $P(m)$ as well.
However, it is still possible to filter out the noise from the sequence $P(m)$ and extract the correct result thanks to a connection between Prony's method and Padé approximation.
In fact, one can prove that the exponential terms in $P(m)$, using the Z-transform, correspond to the poles in the Padé approximation for $f(z) := \sum_{m=0}^{M-1} P(m) z^{-m}$.
From \cite{Cu.Ts.ea:fai:17} we learn that modeling the Padé approximant with additional poles helps to model out the added noise.
The extra poles model the noise and push the others closer to their true locations \cite{Cu.Le:spa:16}.
This requires the use of additional samples of the sequence $P(m)$.
We thus approximate the sequence $P(m)$ with more exponentials than needed, but the extra terms serve to model the noise.

Depending on the amount of noise present in the signal $\boldsymbol{x}$, we have to consider a larger $M$ thus computing more shifted Fourier coefficients ${_u^s}\boldsymbol{X}$.
The number of samples in the sequence $P(m)$ depends on the amount of noise and the number of colliding frequencies we have to discern.
From the sequence $P(m)$ we extract the most prominent exponentials and resolve the aliasing for them.
An example is given in Section \ref{numerical}.

%----------------------------------------------------%
% NUMERICAL
%----------------------------------------------------%
\section{Numerical experiments}\label{numerical}

In this section we apply the proposed technique to some synthetic signals.
At first we focus on the aspect of frequencies colliding because of the undersampling factor $u$.
We consider the following model  

\begin{equation}
  x_l = \sum_{i=1}^{N} \alpha_i e^{ \frac{2 \pi \mathrm{i}}{R} \mu_i l} \quad l=0,\ldots,L-1, 
  \label{examplesignal}
\end{equation}

and we generate three different signals with parameters reported in Table \ref{tableparameters}.

\renewcommand{\arraystretch}{1.3}
\setlength{\tabcolsep}{8pt}
\begin{table}[ht]
  \caption{\it Parameters of the $3$ generated signals.}
  \label{tableparameters}
  \centering
  \begin{tabular}{|c|c|c|c|}
    % \hline
    \multicolumn{4}{c}{Signal $1$} \\
    \hline
    $N$ & $R$ & $\mu_i$  & $\alpha_i$ \\
    \hline 
    $1$ & $1$kHz & $125$ & $1$ \\
    \hline
    \multicolumn{4}{c}{Signal $2$} \\
    \hline
    $N$ & $R$ & $\mu_i$  & $\alpha_i$ \\
    \hline 
    $2$ & $1$kHz & $125$ \quad $165$ & $1$ \quad $e^{\mathrm{i} \pi / 3}$ \\
    \hline
    \multicolumn{4}{c}{Signal $3$} \\
    \hline
    $N$ & $R$ & $\mu_i$  & $\alpha_i$ \\
    \hline 
    $3$ & $1$kHz & $125$ \quad $165$ \quad $245$ & $1$ \quad $e^{\mathrm{i} \pi / 3}$ \quad $e^{\mathrm{i} \pi / 4}$ \\
    \hline
  \end{tabular}
\end{table}

We set $L = 1000$ and perturb the signal with circular white gaussian noise of SNR $= 30$.
If we analyze each signal using the DFT with an undersampling factor $u=50$, the sampling rate becomes $1000/50=20$ Hz and all frequencies collide over the same Fourier bin that represents $5$ Hz.
These collisions caused by the aliasing can be resolved using our method.

\begin{figure}[ht]
  \centering
  \begin{minipage}{0.45\textwidth}
    \centering
    \includegraphics[scale=0.11]{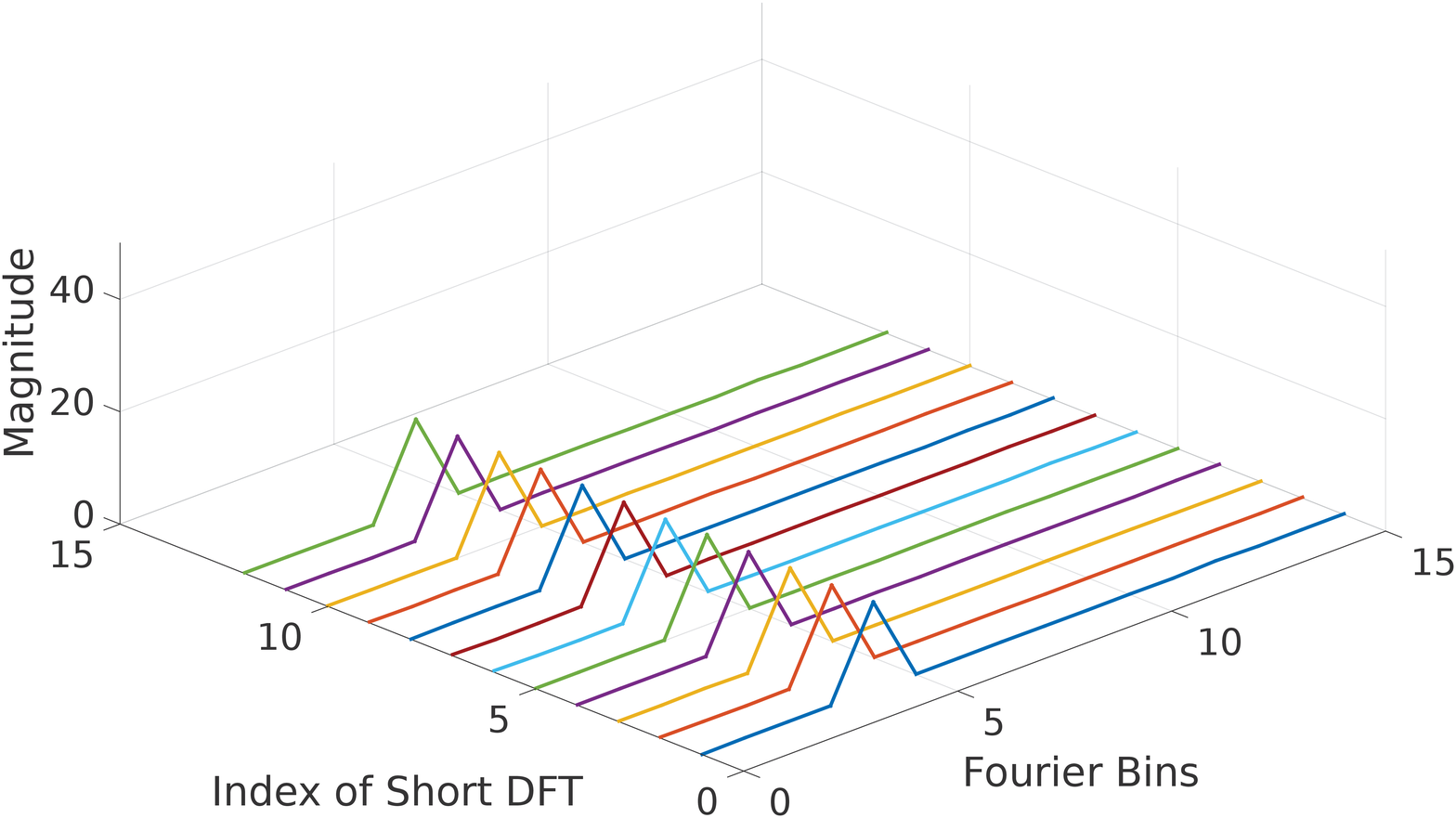}
  \end{minipage}
  \begin{minipage}{0.45\textwidth}
    \centering
    \includegraphics[scale=0.11]{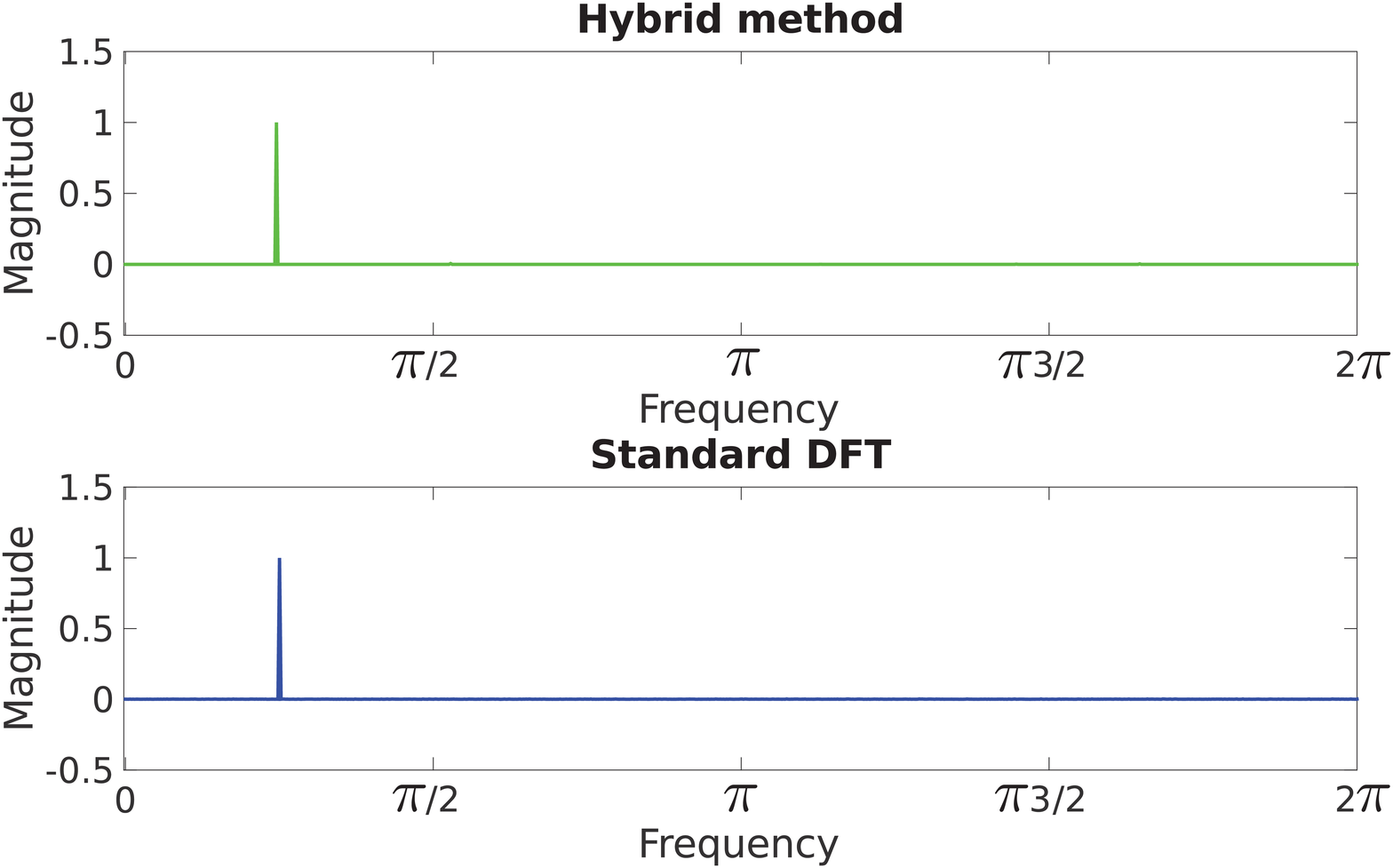}
  \end{minipage}

  \begin{minipage}{0.45\textwidth}
    \centering
    \includegraphics[scale=0.11]{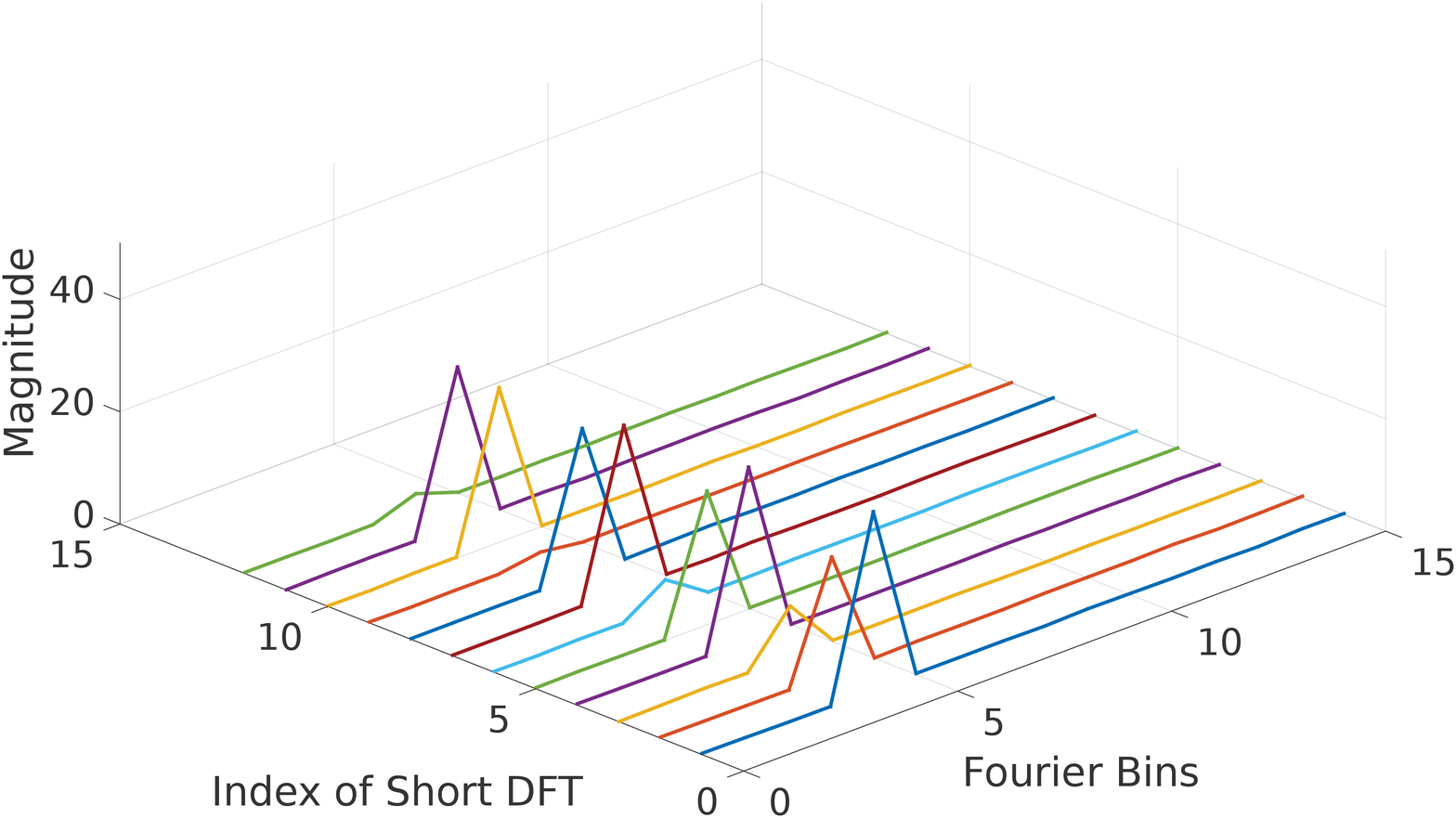}
  \end{minipage}
  \begin{minipage}{0.45\textwidth}
    \centering
    \includegraphics[scale=0.11]{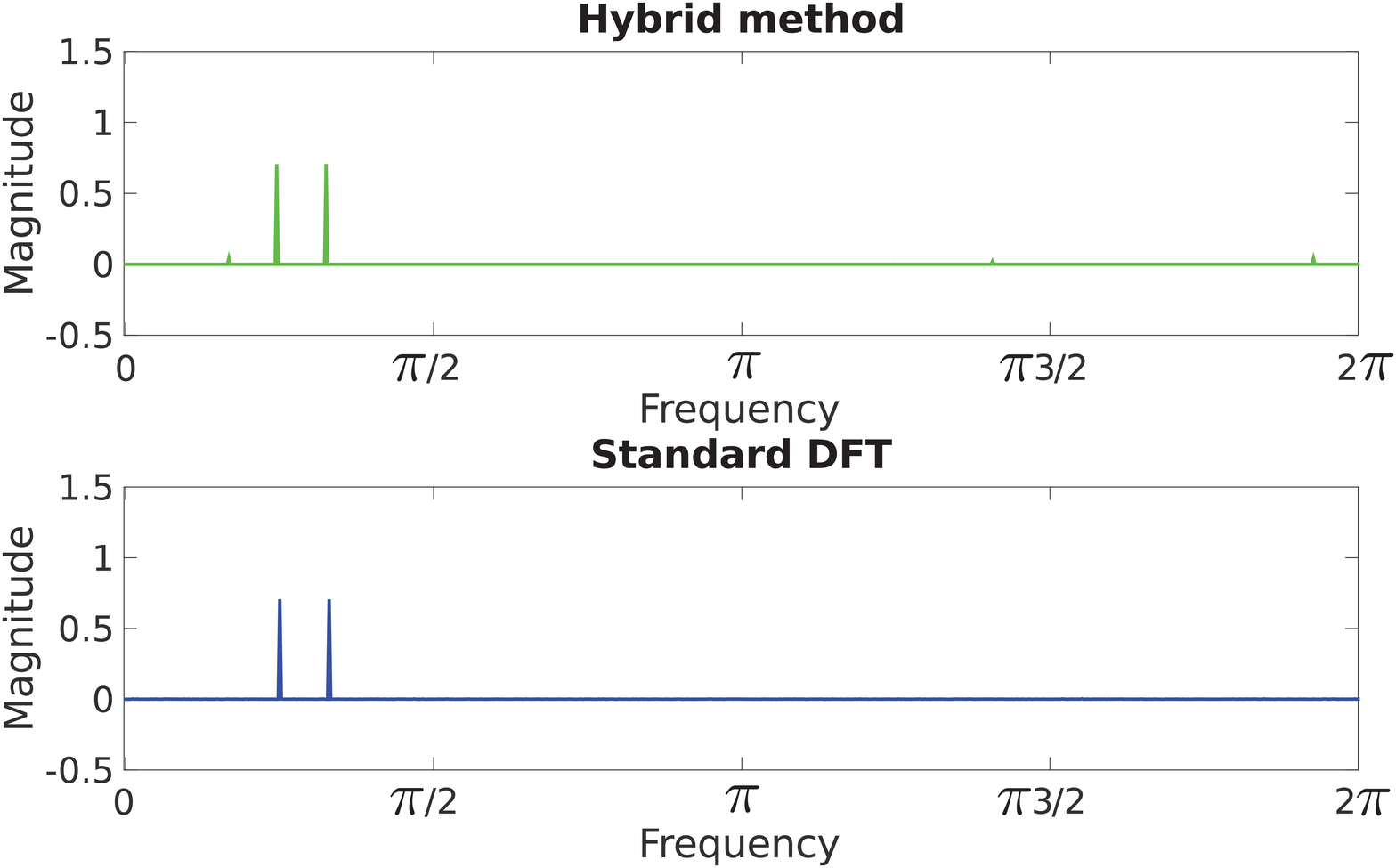}
  \end{minipage}

  \begin{minipage}{0.45\textwidth}
    \centering
    \includegraphics[scale=0.11]{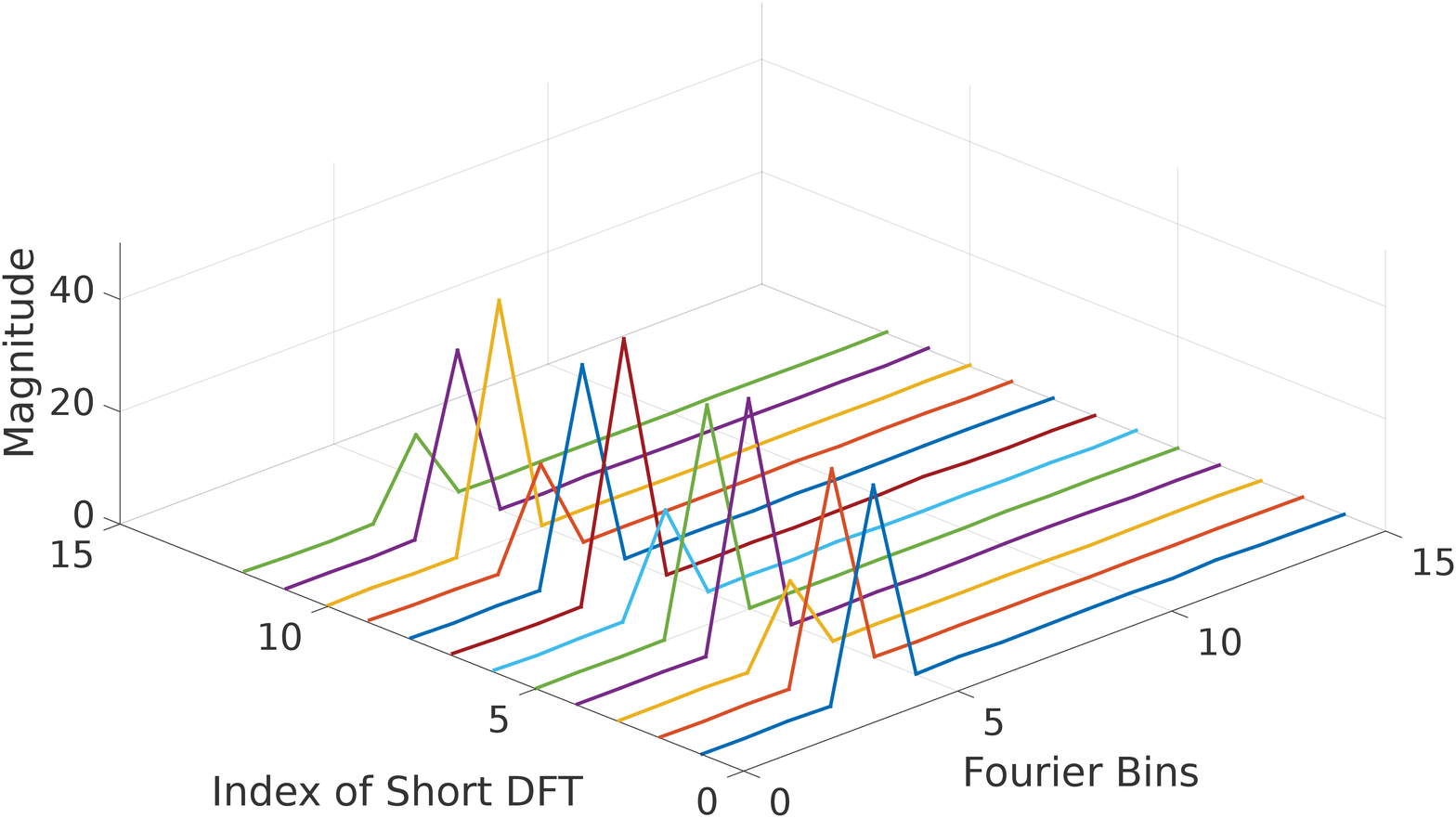}
  \end{minipage}
  \begin{minipage}{0.45\textwidth}
    \centering
    \includegraphics[scale=0.11]{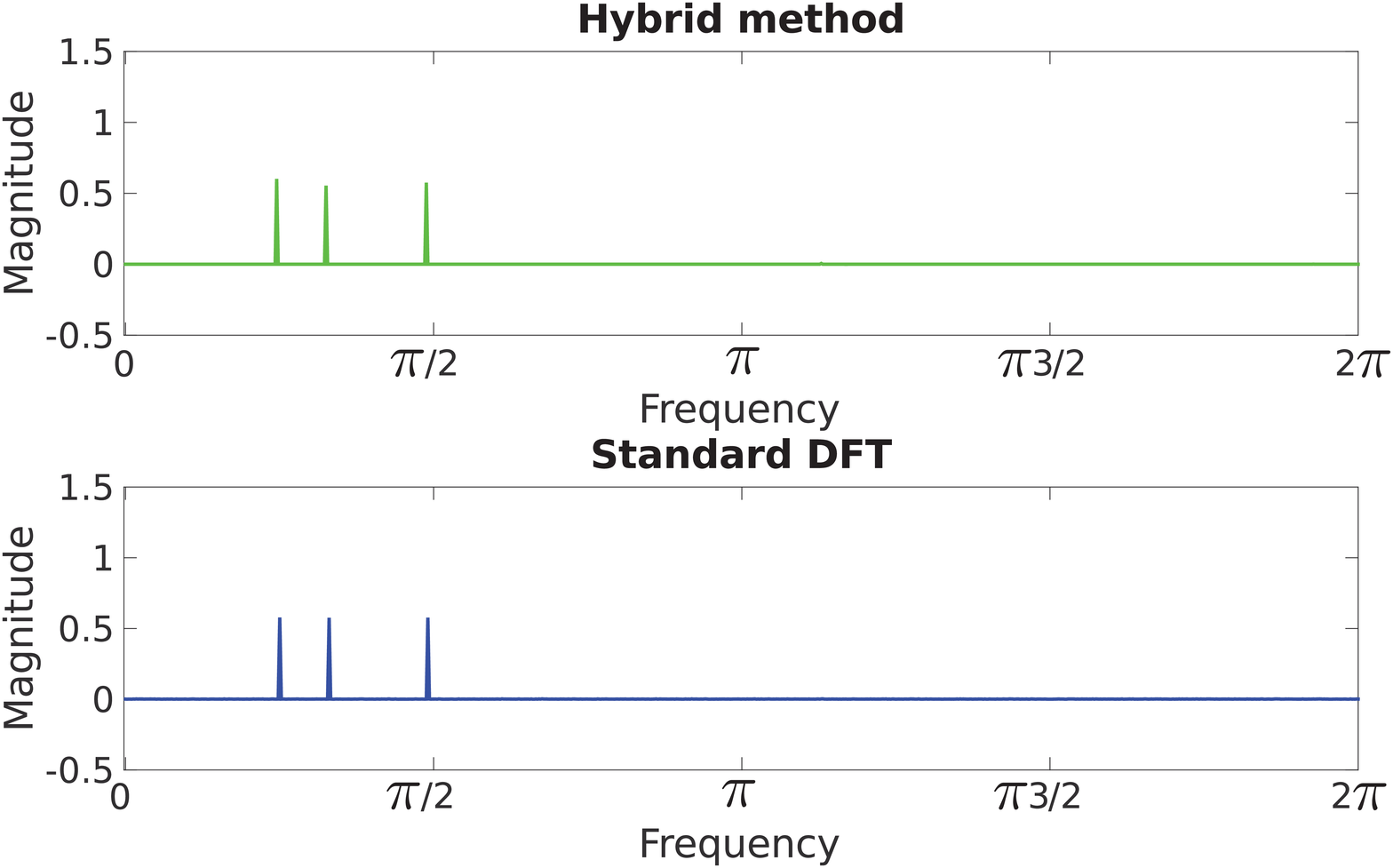}
  \end{minipage}
  \caption{ \textit{ The three generated signals with parameters reported in Table \ref{tableparameters}. 
  On the left hand side the $\lvert{_u^{ms}}\boldsymbol{X}\rvert$ for $m = 0,\ldots,11$ are shown.
  The right hand side shows the results of both the standard DFT and the proposed hybrid method.
  }
  }
  \label{shortDFTs}
\end{figure}

\begin{figure}[ht]
  \centering
  \includegraphics[scale=0.30]{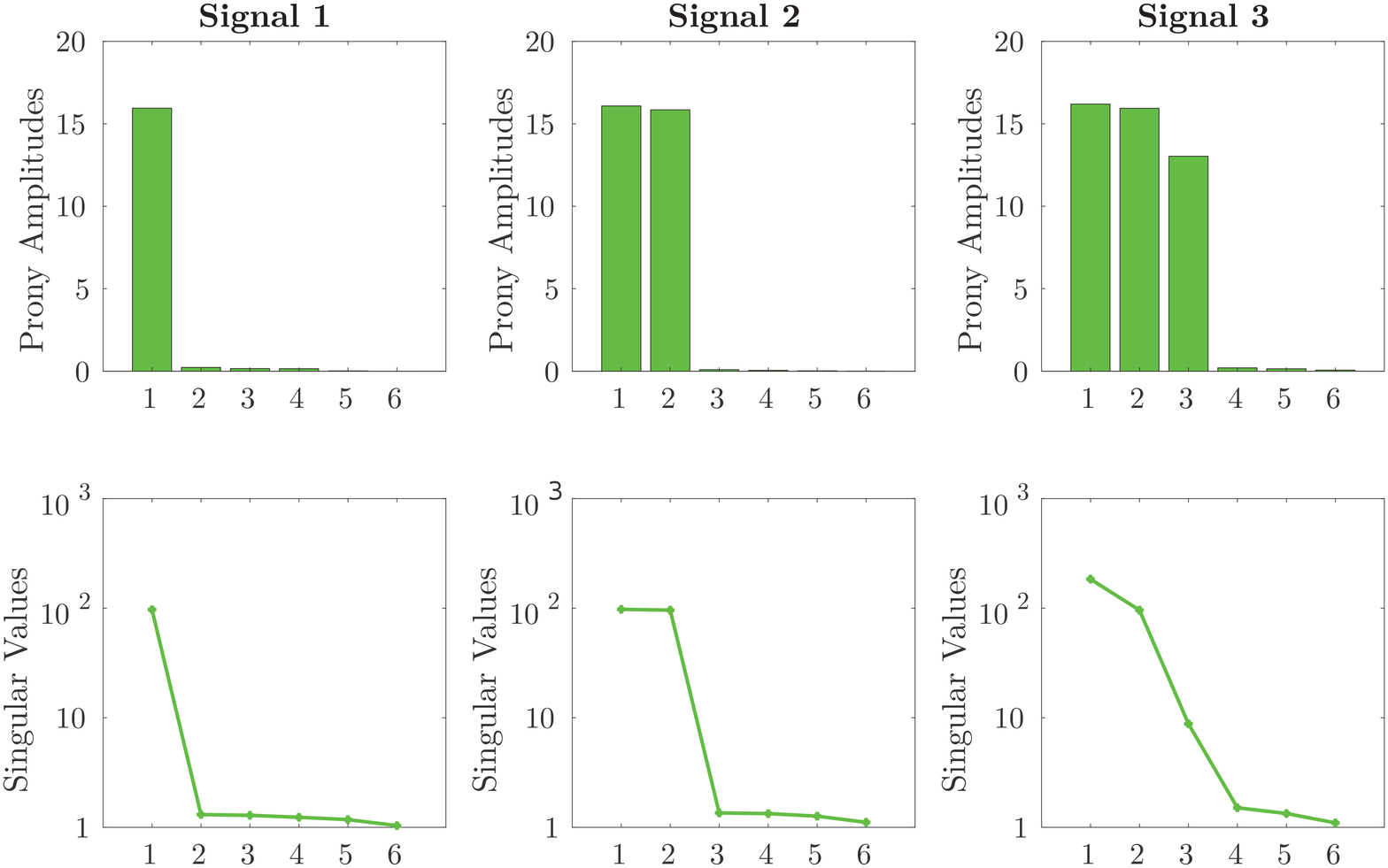}
  \caption{\textit{The first row shows the amplitudes retrieved from the Prony sequences $P(m)$ for each of the three analyzed signals at the $5$ Hz colliding frequency.
  The second row shows the number of meaningful components by reporting the singular values of the Hankel matrices created from the same sequences $P(m)$. }}
  \label{pronyamplitudes}
\end{figure}

We set a shift factor $s=17$ and $M = 12$. 
Due to the limited sampling time in \eqref{examplesignal} because $L=1000$, each undersampled DFT is computed using a total of $\lfloor (L - (s-1)*M ) / u \rfloor = 16$ samples.
For each of the three signals the respective ${^{ms}_u}\boldsymbol{X}, m = 0,\ldots,11$ are portrayed on the left hand side in Figure \ref{shortDFTs}.
On the right hand side of Figure \ref{shortDFTs} the results of the proposed method and the standard DFT using $1000$ samples are reported.
Figure \ref{pronyamplitudes} shows the detail of the sequence $P(m)$ for the frequencies colliding at $5$Hz for each of the signals.

\begin{figure}[ht]
  \centering
  \includegraphics[scale=0.17]{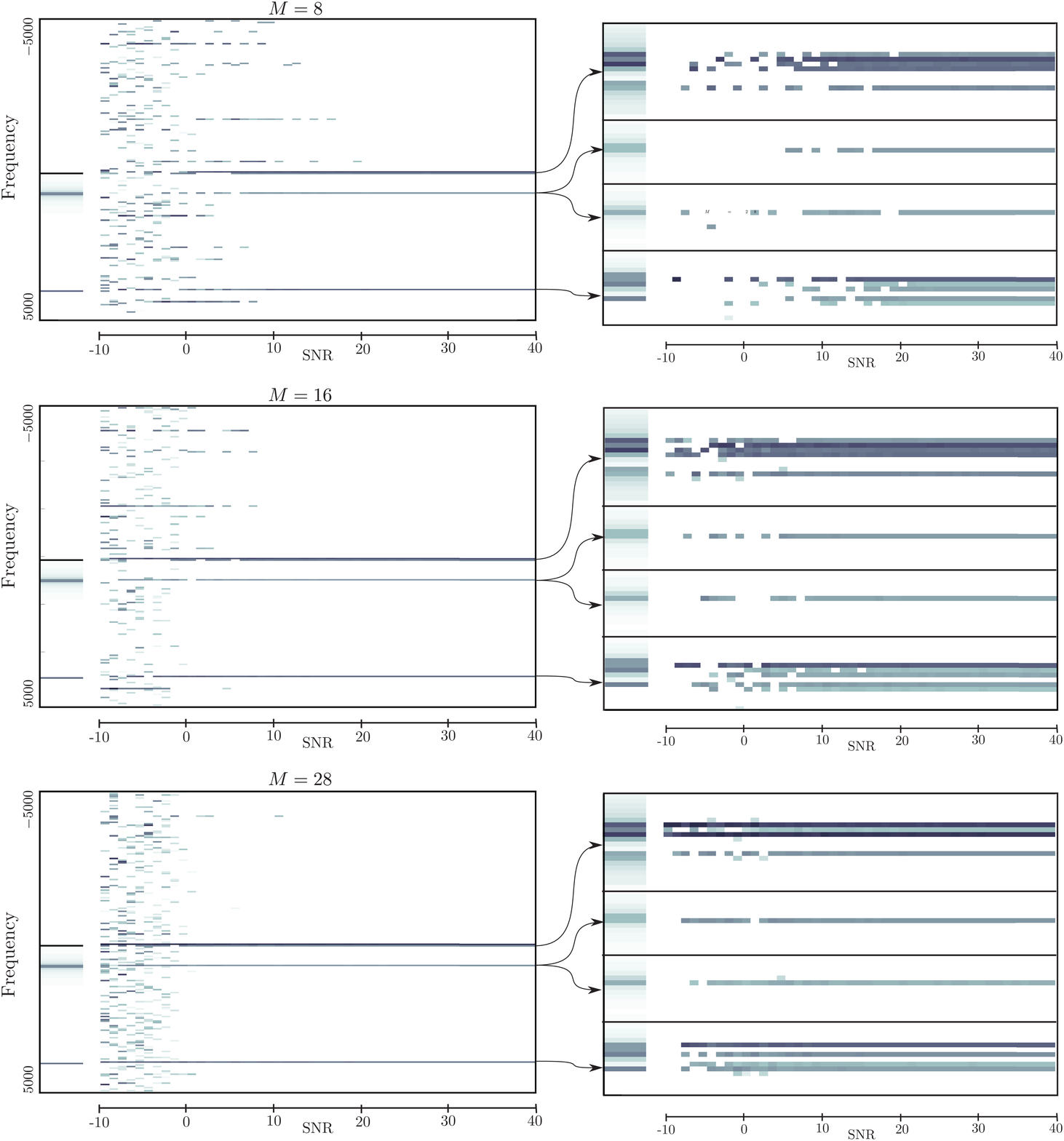}
  \caption{\textit{From top to bottom the same signal perturbed with different noise levels is analyzed using different values of $M$.
  On the left hand side of each figure the output of a noiseless DFT is reported for comparison.
  On the right column are reported detailed regions around the $8$ frequencies of the parameters $\mu_i$ reported in Table \ref{tableparameters2}. }}
  \label{examplelongDFT}
\end{figure}

As second example we consider a signal following model \eqref{examplesignal} with parameters $\mu_i$ and $R$ reported in Table \ref{tableparameters2} and random $\alpha_i$ with $ 0.5 \leq \| \alpha_i \| \leq 1.5$, $N=8$ and $L = 2^{16}$.
The standard DFT computed using all the consecutive samples has resolution $0.1526$Hz, enough to discern the frequency content.
We set $u = 142$, $s = 7 $ and we perform several analyses for different values of $M = 8, 16, 28$, and increasing noise level.
When $M = 28$ the method uses a total of only $ \lfloor (L - (s-1)*M ) / u \rfloor \cdot M = 12824$ samples but maintains a frequency resolution of  $10000 / ( \lfloor ( L - (s-1) \cdot M ) / u  \rfloor \cdot u) = 0.1538$Hz. On the contrary, the standard DFT using the same amount of consecutive samples has a frequency resolution of $10000/12824 = 0.7738$Hz which makes close frequencies indistinguishable. 
Figure \ref{examplelongDFT} illustrates the results for different values of $M$ and increasing noise level up to SNR $= -10$.
The regions around the $\mu_i$ are enlarged in the right hand side of Figure \ref{examplelongDFT}.

\renewcommand{\arraystretch}{1.3}
\setlength{\tabcolsep}{8pt}
\begin{table}[ht]
  \caption{\it Parameters of the second example.}
  \label{tableparameters2}
  \centering
  \begin{tabular}{|c|cccccccc|}
    \hline
    $ R $ & \multicolumn{8}{c|}{ 10000}  \\
    \hline 
    $ \mu_i $ & 100 & 100.3 & 100.92 & 4000 & 4000.3 & 4000.7 & 765 & 787 \\
    \hline
  \end{tabular}
\end{table}

%----------------------------------------------------%
% CONCLUSION
%----------------------------------------------------%
\section{Conclusions}\label{conclusions}

In this paper we propose a novel technique in the field of sparse Fourier methods.
The proposed method uses a sampling scheme that collects samples from several shifted undersampled versions of the original signal.
In this way it can also be seen as a super-resolution technique due to the fact that it achieves a high resolution from sets of samples collected at a lower sampling rate.
The different sample streams could also come from different devices.

The high frequency content is retrieved thanks to a Prony related technique which interpolates the Fourier coefficients affected by aliasing.
The matching of the Fourier and Prony solutions returns the original frequency content.

Using fewer samples than the original DFT scheme, it is still possible to obtain a final frequency resolution comparable to the original one.
The result is a highly parallelizable sparse algorithm that allows to collect data from parallel acquisition systems with a low sampling rate.

%----------------------------------------------------%
% BIBLIOGRAPHY
%----------------------------------------------------%
\clearpage
\nocite{*}
\bibliographystyle {plain}
\bibliography{mylib}

\end{document}